# Incremental Model Building Homotopy Approach for Solving Exact AC-Constrained Optimal Power Flow

Amritanshu Pandey
ECE, Carnegie Mellon University
amritanp@andrew.cmu.edu

Aayushya Agarwal
ECE, Carnegie Mellon University
aayushya@andrew.cmu.edu

Larry Pileggi
ECE, Carnegie Mellon University
pileggi@andrew.cmu.edu

**Abstract**

*Alternating-Current Optimal Power Flow (AC-OPF) is framed as a NP-hard non-convex optimization problem that solves for the most economical dispatch of grid generation given the AC-network and device constraints. Although there are no standard methodologies for obtaining the global optimum for the problem, there is considerable interest from planning and operational engineers in finding a local optimum. Nonetheless, solving for the local optima of a large AC-OPF problem is challenging and time-intensive, as none of the leading non-linear optimization toolboxes can provide any timely guarantees of convergence. To provide robust local convergence for large complex systems, we introduce a homotopy-based approach that solves a sequence of primal-dual interior point problems. We utilize the physics of the grid to develop the proposed homotopy method and demonstrate the efficacy of this approach on U.S. Eastern Interconnection sized test networks.*

## 1. Introduction

Alternating Current Optimal Power Flow (AC-OPF) attempts to solve for the most economical dispatch of the grid while satisfying AC network constraints and enforcing device limits. Posed as a non-convex, non-linear optimization problem, AC-OPF is **NP-hard** with no methodology for guaranteeing a global optimum in polynomial time for general networks even though some work [1]-[2] have found the global optimum for a restricted set of small networks.

Instead, many methodologies have been proposed and developed to solve for the local optimum or an approximation of the AC-OPF problem due to increasing demand from grid planners and operators [3]. These methodologies can be divided into three main classes: i) those that solve for the local optimum of the original AC-OPF problem [4]-[5]; ii) those that solve the linear approximation of the AC-OPF problem (e.g. Decoupled-OPF) [6]-[8]; and iii) those that solve the convex relaxation or restriction of the AC-OPF problem [9]-[12].

Of these three classes of methodologies, in this paper, we focus on obtaining the local optimum solutions for the AC-OPF problems primarily because the local optimum ensures a feasible dispatch with AC network constraints satisfied and device limits enforced. In planning and operation, security and reliability are considered more critical than the cost of the electricity and obtaining a sub-optimal feasible solution to AC-OPF is generally preferred to an AC-infeasible solution (due to relaxations) that is more economical.

Due to these factors, there has been a recent surge in methods that locate local optimum solutions for the AC-OPF problem [4]-[5]. In particular, a large-scale effort was driven by an Advanced Research Project Agency – Energy (ARPA-E) grid optimization challenge [13]. The goal of the competition was to locate the most economical dispatch of the AC-OPF while ensuring feasibility across a set of contingencies (also known as Security Constrained Optimal Power Flow problem). Many state-of-the-art methods for obtaining local optimum solutions for non-convex optimization problems were presented in the competition [13], with the most successful ones based on standard nonlinear optimization toolboxes such as the Matpower Interior Point Solver (MIPS) [14] and Interior Point Optimizer (IPOPT) [15]. The approaches utilized power-system specific knowledge to achieve robust convergence. An alternative approach made use of a circuit-based formulation with currents and voltages as state variables [5] that applied circuit-simulation methods to achieve robust convergence. Almost all of these approaches demonstrated the ability to solve real-size AC-OPF problems within the SC-OPF problem; nonetheless, the results of the challenge [13] suggest that there is value in further exploration. Most importantly, since the optimization framework also serves as a basis for planning engineers, a robust approach cannot always rely on starting with good initial conditions.

Various previous works have attempted to develop a robust AC-OPF approach by utilizing a homotopy method technique [16] that traces a path in the solution space by incrementally modifying the non-linear constraints from a trivial problem to the original AC-OPF. This method requires a continuous path to the solution with a valid solution at each increment. However, previous works [18]-[24] are unable to ensure a valid feasible solution at each increment and few create discontinuities in the homotopy path due to the use of DC-OPF as an initial solution [19]. As a result, many homotopy methods used for the power grid today [18]-[24] are unable to scale to large networks.

One homotopy method that is effective for large-scale power-flow problems is Tx-stepping [25]. Tx-stepping is based on an approach that is analogous to Gmin stepping for circuits [27] and aims to short the series elements (e.g. lines, transformer, etc.) in the system by adding large parallel admittances at first to obtain a trivial solution. Throughout the homotopy path then, the added admittances are gradually removed to solve the original network. Here for grid optimization we employ a more comprehensive incremental formation of the final system starting from a trivial one that not only scales the transmission line and transformer parameters, but also shunts, generators, and loads and corresponding bounds. This effectively corresponds to an incremental formulation of the entire grid from a trivial one that mimics *turning on the grid*. Importantly, this incremental stepping must be performed such that there is a valid grid model at each step, which we ensure by using a slack injection-based technique that maintains feasibility throughout the homotopy method. We refer to this homotopy as **IMB, Incremental Model Building**.

We introduce a framework that implements the IMB approach within the Simulation of Unified Grid Analyses and Renewables (SUGAR) engine [5], [17]. We overcome convergence challenges that general-purpose optimization toolboxes struggle with by utilizing the known properties of the physics-based grid models to ensure convergence, even when confronted with ill-suited initial conditions. Importantly, the proposed homotopy approach is generic and can be applied to achieve robust convergence for other grid optimization problems as well.

Section 3 covers the AC-OPF formulation along with primal-dual interior point (PDIP) solution methodology. Section 4 introduces the novel homotopy method and Section 5 describes the algorithmic methodology following the proposed approach. We conclude the paper with AC-OPF results for large and complex networks and compares it against state-of-the-art standard optimization toolboxes in Section 6.

## 2. Nomenclature

| | |
|---|---|
| $a_g, b_g, c_g$ | coefficients of quadratic cost function for generator $g$. |
| $I_{G,i}^R, I_{G,i}^I$ | Real and imaginary currents injections by a generator. |
| $I_{D,i}^R, I_{D,i}^I$ | Real and imaginary currents injections by a load. |
| $V_i^R, V_i^I$ | Real and imaginary voltage at node $i$. |
| $V_{ik}^R, V_{ik}^I$ | Real and imaginary voltage across nodes $i$ and $k$. |
| $G_{ik}^Y, B_{ik}^Y$ | Components of admittance matrix. |
| $P_g^G, Q_g^G$ | Real and reactive power of generation $g$. |
| $P_d^D, Q_d^D$ | Real and reactive power of constant power demand $d$. |
| $V_i^{SQ}$ | Square of the voltage magnitude at bus $i$. |
| $I_e^{ik}$ | Current flow in a series element $e$. |
| $S_e^{ik}$ | Apparent power flow in a series element $e$. |
| $\underline{P_g}, \overline{P_g}$ | Bounds on generator real power. |
| $\underline{Q_g}, \overline{Q_g}$ | Bounds on generator reactive power. |
| $\underline{V_i}, \overline{V_i}$ | Bounds on voltage magnitude at bus $i$. |
| $t_{r,i}, \theta_i$ | Turns ratio and phase shift angle for transformer $i$. |
| $\overline{I_e^{ik}}$ | Bound on current flow in a series element $e$ between node $i$ and $k$. |
| $\overline{S_e^{ik}}$ | Bound on apparent power flow in a series element $e$ between node $i$ and $k$. |
| $F(x)$ | Original system of non-linear equations. |
| $\mathcal{H}(x,v)$ | Homotopy parameterized system of non-linear equations. |
| $G(x)$ | System of non-linear equations with a trivial solution. |
| $v$ | Homotopy factor. |
| $\gamma$ | Scaling factor of admittance during homotopy. |
| $c(v)$ | Homotopy path as a function of homotopy factor $v$. |
| $\overline{\mu}, \underline{\mu}$ | Vector of slack variables corresponding to the upper and lower bound of states vector $x$, respectively. |
| $\lambda$ | Vector of dual variables for equality constraints. |
| $x, \overline{x}, \underline{x}$ | System state vector along with the vector of upper and lower bounds, respectively. |

## 3. Alternating Current OPF Formulation

### 3.1. Current Voltage Formulation

Traditionally, the non-linear AC-network constraints of an AC-OPF problem are represented by power mismatch at each bus; however, the trigonometric terms that model network line flows and network constraints introduces steep non-linearities in the solution space when scaling to large systems. A growing number of approaches [3], [5] are adopting a current-voltage formulation in which network constraints are expressed as Kirchhoff's Current Laws (KCL). The current-voltage (I-V) models have better convergence properties for the AC-OPF problem than the traditional power mismatch-based formulations due to the reduced non-linearities in the network constraints [3]. Unlike the traditional formulation where the network constraints are non-linear with trigonometric terms, the I-V formulation has linear network constraints with non-linearities stemming from the injection models (loads and generation). Therefore, to minimize non-linearities in the formulation, we will utilize the I-V representation of the AC-OPF problem.

The optimization formulation for the AC-OPF problem of a power grid network $\mathcal{N}$ consisting of a set of generators $\mathcal{G}$ and load demands $\mathcal{D}$ connected to a set of buses or nodes $B$ in the grid is given by (1). The nodes $B$ in the system are connected by a set of network elements, $\mathcal{T}_X$ and $X_F$. The objective function for AC-OPF of minimizing generation cost is given by $\mathcal{F}_c(P^G)$, defined by sum of quadratic functions given by a set of coefficients $\{a, b, c\}$. The equality constraints represent the network constraints given by (1b)-(1h). Inequality constraints represent physical bounds on the devices given by (1i)-(1l). Some bounds include system limits that are based on grid stability; for example flow constraints can be both thermal or stability based (see 1l) and are further discussed in Section 4.4. In the problem formulation, system topology parameters, cost, and device and voltage bounds are given parameters whereas generator real and reactive power output along with voltage-setpoint are the decision variables ($x_d \subseteq x = \{P^G, Q^G, V_G^{SQ}\}$) over whose range the problem is optimized. More detailed formulations can also include transformer tap and phase-shifters ($t_r, \theta$) along with shunt positions as decision variables.

$$\min_{P^G} \mathcal{F}_c(P^G) = \sum_{g=1}^{|\mathcal{G}|} [a_g + b_g P_g^G + c_g (P_g^G)^2] \quad (1a)$$

subject to:

$$I_{G,i}^R - I_{D,i}^R = \mathbf{Re}\left\{\sum_{k=1}^{|B|}(V_{ik}^R + jV_{ik}^I)(G_{ik}^Y + jB_{ik}^Y)\right\} \quad \forall i \in B \quad (1b)$$

$$I_{G,i}^I - I_{D,i}^I = \mathbf{Im}\left\{\sum_{k=1}^{|B|}(V_{ik}^R + jV_{ik}^I)(G_{ik}^Y + jB_{ik}^Y)\right\} \quad \forall i \in B \quad (1c)$$

$$I_{G,i}^R = \sum_{g=1}^{|\mathcal{G}(i)|} \frac{-P_g^G V_i^R + Q_g^G V_i^I}{V_{R,i}^2 + V_{I,i}^2} \quad \forall i \in B \quad (1d)$$

$$I_{G,i}^I = \sum_{g=1}^{|\mathcal{G}(i)|} \frac{-P_g^G V_i^I - Q_g^G V_i^R}{V_{R,i}^2 + V_{I,i}^2} \quad \forall i \in B \quad (1e)$$

$$I_{D,i}^R = \sum_{d=1}^{|\mathcal{D}(i)|} \frac{P_d^D V_i^R + Q_d^D V_i^I}{V_{R,i}^2 + V_{I,i}^2} \quad \forall i \in B \quad (1f)$$

$$I_{D,i}^I = \sum_{d=1}^{|\mathcal{D}(i)|} \frac{P_d^D V_i^I - Q_d^D V_i^R}{V_{R,i}^2 + V_{I,i}^2} \quad \forall i \in B \quad (1g)$$

$$V_i^{SQ} = V_{R,i}^2 + V_{I,i}^2 \quad \forall i \in B \quad (1h)$$

$$\underline{P_g} \leq P_g^G \leq \overline{P_g} \quad \forall g \in \mathcal{G} \quad (1i)$$

$$\underline{Q_g} \leq Q_g^G \leq \overline{Q_g} \quad \forall g \in \mathcal{G} \quad (1j)$$

$$\left(\underline{V_i}\right)^2 \leq V_i^{SQ} \leq \left(\overline{V_i}\right)^2 \quad \forall i \in \mathcal{N} \quad (1k)$$

$$\left(I_e^{ik}\right)^2 \leq \left(\overline{I_e^{ik}}\right)^2 \quad \forall e \in \{\mathcal{T}_X, X_F\} \quad (1l)$$

With the non-linear optimization formulation devised, we can construct a framework for solving for a local minimum using Primal-Dual Interior Point method.

### 3.2. Solution Framework for AC-OPF Problem

To solve the non-convex optimization problem given in (1), we implement the primal-dual interior point (PDIP) approach (see [26] for details). The PDIP approach is the basis for further heuristics developed within this paper that enables robust and scalable convergence.

PDIP algorithms apply the search direction from Newton's method to iteratively solve a set perturbed (Karush Kuhn Tucker) KKT conditions. The perturbed KKT conditions represent a relaxation of first-order optimality conditions that are necessary to obtain a local optimum for the AC-OPF problem. We begin by formulating the Lagrangian:

$$\mathcal{L}(x, \lambda, \mu) = \mathcal{F}_c(x) + \lambda^T g(x) + \mu^T h(x) \quad (2)$$

where, $g(x)$ represents the vector of equality constraints given by (1b)-(1h) (in case of AC-OPF) and $h(x)$ represents the vector of inequality constraints (1i)-(1l). $x$ represents the vector of primary variables whereas $\lambda$ and $\mu$ represent the vector of dual and slack variables corresponding to equality and inequality constraints, respectively. In the case of AC-OPF, $x$ is a vector of grid states that include real and reactive power output of the generators ($P^G, Q^G$), vector of real and imaginary components of grid voltages ($V_R, V_I$), vector of line and transformer current and power flows ($I^{ik}, S^{ik}$) and a vector of generator set-points ($V_G^{SQ}$). The set of perturbed KKT conditions (perturbed due to the relaxation in complementarity slackness constraint) corresponding to (2) can now be written as:

$$\nabla_\lambda \mathcal{L} = g(x) = 0 \quad (3a)$$

$$\nabla_x \mathcal{L} = \nabla_x^T g(x)\lambda = 0 \quad (3b)$$

$$\overline{\mu} \circ (x - \overline{x}) + \epsilon = 0 \quad (3c)$$

$$\underline{\mu} \circ (x - \underline{x}) - \epsilon = 0 \quad (3d)$$

$$\mu \geq 0 \quad (3e)$$

$$\underline{x} \leq x \leq \overline{x} \quad (3f)$$

where $\circ$ represents element-wise multiplication of the vector elements and $\epsilon$ represents a vector of complementarity slackness tolerance.

To obtain a stationary point ($x^*, \lambda^*, \mu^*$) that satisfies (3), PDIP linearizes and iteratively solves the equations corresponding to (3a)-(3d) using Newton's method. Between two NR iterates, primal and dual feasibility given by (3e)-(3f) is satisfied through heuristics. The linearized matrix at iteration $k$ can be written as follows:

$$\begin{bmatrix} H_k & \nabla_x^T g(x^k) & I \\ \nabla_x g(x^k) & 0 & 0 \\ U_k & 0 & X_k \end{bmatrix} \begin{bmatrix} \Delta x \\ \Delta \lambda \\ \Delta \mu \end{bmatrix} = - \begin{bmatrix} \nabla f_c(x^k) + \nabla_x^T g(x^k) \lambda + \mu^T \\ \nabla_x^T g(x^k) \\ \epsilon \end{bmatrix} \quad (4)$$

where $U \coloneqq \text{diag}(\mu)$, $X \coloneqq \text{diag}(x)$ and $H = \nabla_{xx}^2 \mathcal{L}(x^k, \lambda^k, \mu^k)$. In general, the linearized matrix is reduced and $\mu$'s eliminated to solve lower dimension matrix within the inner loop of the Newton's method. A descent direction is guaranteed by ensuring that the top left block (given by $H_k$) of the iteration matrix is positive definite [15] and step-length of the iterate is often controlled via backtracking line search [26].

While this straightforward PDIP solution methodology could be used to solve small networks, it is unlikely to be effective for complex large networks with steep nonlinearities. Therefore, to ensure convergence for such realistic networks, heuristics are needed. Most commercial tools have such heuristics embedded within them. For instance, MATLAB optimization toolbox makes use of interior trust region methods [32], whereas IPOPT makes use of filter line search-based methods along with second order correction terms [15]. Nonetheless, these methods also diverge for some of the hardest networks as these do not intrinsically utilize grid physics within their solvers. Therefore, to ensure robust convergence of hardest AC-OPF networks, we augment the naïve PDIP algorithm with a novel IMB homotopy method that utilizes the physics of the electric grid to robustly obtain a local optimum.

## 4. Incremental Model Building (IMB) Homotopy for AC-OPF

### 4.1. Homotopy Methods

Homotopy methods [16] are commonly used to solve complex large-scale non-convex optimization problems pertaining to many applications. However, so far, the use of such methods for solving the power grid optimization problems is limited [17]-[24]. These methods embed a scalar homotopy factor, $v$, into the problem in an effort to relax the non-linearities. Varying the homotopy factor, $v: 1 \rightarrow 0$, replaces the original problem with a set of sub-problems represented by each increment of the homotopy factor and are sequentially solved with the following properties: i) the first sub-problem has a trivial solution and ii) each subsequent problem has a solution very close to the solution of the prior sub-problem wherein the final subproblem represents the original problem. This second property exploits the Newton-Raphson quadratic convergence properties thereby allowing faster convergence. Mathematically this can be described via the following expression:

$$\mathcal{H}(x, v) = (1 - v)\mathcal{F}(x) + v\mathcal{G}(x) \quad (5)$$

where $v \in [0, 1]$.

The method begins by embedding the homotopy factor into the original problem, effectively replacing the original problem $\mathcal{F}(x) = 0$ with $\mathcal{H}(x, v) = 0$. The equation set $\mathcal{G}(x)$ is a representation of the system that has a trivial solution. The homotopy factor $v$ has the value of 1 for the first sub-problem that corresponds to a trivial problem $\mathcal{G}(x)$. Iteratively the homotopy factor, $v$, is reduced to 0 which then represents the original problem $\mathcal{F}(x)$. In the homotopy path between the trivial problem $\mathcal{G}(x)$ and the original problem $\mathcal{F}(x)$ lies a sequence of sub-problems that traces a path in the solution space. Importantly, the homotopy method requires a valid solution at each step of the homotopy path to ensure convergence.

Previously, the homotopy method has been proven to be an effective tool to achieve robust convergence for complex power flow analysis problems in the form of a Tx-stepping in [25]. However, power flow analysis assumes a fixed generation and solely applying Tx-stepping in an AC-OPF setting proves ineffective as generation is variable in the problem definition. In this paper, we extend the Tx-stepping homotopy method to develop a broader more comprehensive homotopy method to solve optimization problems. Specifically, we target robust solution of AC-OPF problems using our Incremental Model Building (IMB) approach.

### 4.2. Tx-Stepping Homotopy Method

In previous work [25] we developed a homotopy method, **Tx stepping,** that solved for any positive sequence or three-phase power flow (PF) problem independent of the given initial conditions. Tx-stepping is based on embedding a homotopy factor in the linear transmission network to virtually short the grid. This homotopy approach begins by solving an almost shorted system and then it gradually decreases the homotopy factor to obtain the solution to the original problem. We briefly describe how grid models are modified based on Tx-stepping method.

**4.2.1. Transmission line models.** In Tx-stepping, the series elements in the system (transmission lines, transformers etc.) are "virtually" shorted by adding a large conductance ($v\gamma G_i$) and a large susceptance ($v\gamma B_i$) in parallel to each transmission line and transformer model $i$ in the system parameterized by the homotopy factor:

$$\forall i \in \{\mathcal{T}_X, X_F\} : \hat{G}_i = G_i + v\gamma G_i \quad (6)$$

$$\forall i \in \{\mathcal{T}_X, X_F\} : \hat{B}_i = B_i + v\gamma B_i \quad (7)$$

where, $X_F$ is the set of all transformers and $\mathcal{T}_X$ is the set of all the transmission lines in the system. $G$ and $B$ are the original line and transformer impedances and the $\hat{G}$ and $\hat{B}$ are the parameterized admittances during Tx-Stepping that are used while iterating from trivial problem to the original problem. The parameter $\gamma$ is used as a scaling factor for the conductances ($G$) and susceptances ($B$). Evidently, the first homotopy factor will increase the conductance and susceptance of all the branches and transformers $((\hat{G}_i \text{ and } \hat{B}_i), \forall i \in \{\mathcal{T}_X, X_F\})$ thereby almost shorting the system and incrementally the homotopy factor will decrease these values to the original network state at $v = 0$.

**4.2.2. Transformer phase shifters and taps.** A "virtually short" power system (at a homotopy factor of 1) will drive all the voltage magnitude and angle of all the buses to the same value driven by the reference bus. To ensure transformer taps and phase shifter angles are compliant, their turns ratios, $\hat{t}_r$ and phase shift angles, $\hat{\theta}$ correspond to a magnitude of 1 pu and 0°, respectively. Subsequently, the homotopy factor $v$ is reduced to constrain the transformer tap and phase shifters to their original settings. This can be mathematically expressed as follows:

$$\forall i \in X_F : \hat{t}_{r,i} = t_{r,i} + v(1 - t_{r,i}) \quad (8)$$

$$\forall i \in X_F : \hat{\theta}_i = \theta_i - v\theta_i \quad (9)$$

Note that with following changes to the transmission line and transformer models, we have shown to solve any power flow problem robustly, independent of the choice of the initial conditions [17]. However, this approach is insufficient and unlikely to work for AC-OPF as the generation is variable and the physical limits of devices must be enforced. Therefore, to overcome these challenges, we devise a new homotopy method as an augmentation of Tx-stepping method.

## 4.3. Homotopy Models for Optimal Power Flow

Incremental Model Building (IMB) builds upon Tx-stepping to provide a homotopy method that is designed for optimizing power grids. The goal of the IMB homotopy method is to mimic gradually *turning on the grid* from a shorted no-load network. Specifically, we first scale down all the loads in the network to virtually zero in the first homotopy step by embedding a homotopy factor in load and other AC-OPF based grid models, thereby creating a viable solution model that satisfies the physics constraints. From there the embedded homotopy factor iteratively traces a path from the initial sub-problem, depicted in Figure 1 (left) to the final one as shown later in Figure 1 (right). Along with virtually shorting the transmission network as is done in Tx-stepping, the IMB approach embeds the homotopy factor into the load, generator, and shunt models to trace the solution path.

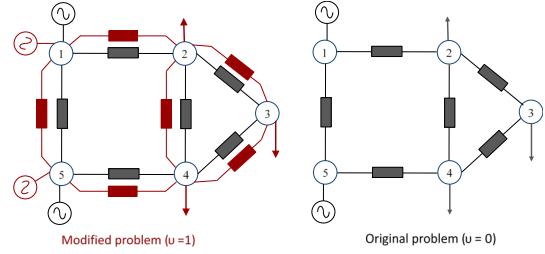

*Figure 1: Schematic of original and modified network during extended Tx-stepping.*

**4.3.1. Load Models.** To mimic turning off the grid at first, we embed a homotopy factor ($v$) that reduces a load's active and reactive powers to virtually zero at first and then gradually scales it back up as it traverses through the homotopy path, as shown by (10)-(11). The initial trivial problem with $\bar{v} = \min(0.0001, v)$, modifies the network to have virtually no load.

$$\hat{I}_{D,i}^R = \sum_{d=1}^{|D(i)|} \frac{(1-\bar{v})(P_d^D V_i^R + Q_d^D V_i^I)}{V_{R,i}^2 + V_{I,i}^2} \quad \forall i \in B \quad (10)$$

$$\hat{I}_{D,i}^I = \sum_{d=1}^{|D(i)|} \frac{(1-\bar{v})(P_d^D V_i^I - Q_d^D V_i^R)}{V_{R,i}^2 + V_{I,i}^2} \quad \forall i \in B \quad (11)$$

**4.3.2. Generator Limits.** During the initial step with the "turned-off grid" homotopy step, the output of generator's real power is very close to zero ($P_g \cong 0, \forall g \in \mathcal{G}$). To ensure that this output is feasible and within the generator real power bounds, we adjust the upper and lower limits of the generator real and reactive power based on the homotopy parameter. This ensures feasible generator operation in the homotopy path even when the original lower bound of the generator real power is much greater than 0 ($\underline{P_g} \gg 0$) or the original higher bound is much lower than 0 ($0 \gg \overline{P_g}$).

$$\overline{\hat{P}_g} = (1-v)\overline{P_g} + v\kappa \quad \forall g \in \mathcal{G} \quad (12)$$

$$\underline{\hat{P}_g} = (1-v)\underline{P_g} - v\kappa \quad \forall g \in \mathcal{G} \quad (13)$$

$$\overline{\hat{Q}_g} = (1-v)\overline{Q_g} + v\kappa \quad \forall g \in \mathcal{G} \quad (14)$$

$$\hat{Q}_g = (1-v)\underline{Q}_g - v\kappa \qquad \forall g \in \mathcal{G} \quad (15)$$

where $\kappa$ is a fixed constant (in our experience a value of 0.5 works well).

### 4.3.3. Shunt models.
Following the IMB technique, shunt devices are initially turned-off and then gradually turned on while traversing through the homotopy path to have its full capacity when solving the original problem:

$$\forall i \in sh: \hat{G}_i^{sh} + j\hat{B}_i^{sh} = (1-v\gamma)(G_i^{sh} + jB_i^{sh}) \quad (16)$$

$(\hat{G}_i^{sh}, \hat{B}_i^{sh})$ are homotopy parameterized shunt impedances that replace the original impedances $(G_i^{sh}, B_i^{sh})$ at bus $i$. Iteratively with the homotopy factor, the shunt admittances are restored to their original value.

## 4.4. Circuit-based Flow Constraints

Power and current-based limits are two common ways in which flow constraints for transmission line and transformers are represented within AC-OPF problem. They either represent thermal or stability limits for the series element. Commonly, current bounds are used for lines and power bounds are used for transformers. We first describe how we implement current bounds. Power bounds can be implemented similarly.

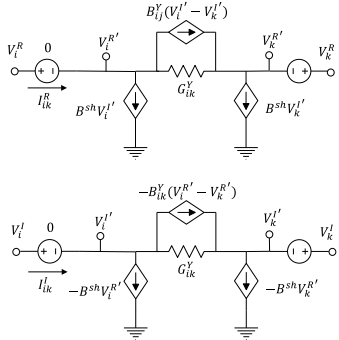

*Figure 2: Circuit representation for measuring transmission line currents.*

Modeling the line and transformer current bounds require computing the real and imaginary currents through any given transmission line or transformer. To calculate the currents as in the case of circuit simulation, we append a voltage source with a voltage of 0 pu (representing an ammeter) in series with the series element, as shown by Figure 2. This results in addition of two unknown states to the set of equations, which corresponds to measured real and imaginary currents ($I_{ik}^R$ and $I_{ik}^I$). Importantly, the voltage measurement source is a linear element that is able to measure the real and imaginary currents, of $I_{ik}^R$ and $I_{ik}^I$ without introducing any non-linearities. Now that these currents are state variables, we can add an additional equation (17) to calculate the square of magnitude of the current through the series element, $(I_e^{ik})^2$, where $\forall e \in \{\mathcal{T}_X, X_F\}$.

$$(I_e^{ik})^2 - \left( \left(I_{ik,e}^R\right)^2 + \left(I_{ik,e}^I\right)^2 \right) = 0 \quad (17)$$

Inserting (17) into the Lagrange function (2) with an associated dual variable, $\lambda_{I,e}$, we can bound the current magnitude to a value $\overline{I_e^{ik}}$ by adding following constraint:

$$\mu_{I,e}\left( (I_e^{ik})^2 - \left(\overline{I_e^{ik}}\right)^2 \right) + \epsilon = 0 \quad (18)$$

Similar methodology can be followed if transformer bounds are given in apparent power units, with an added equation for power calculated from receiving node and sending node:

$$(S_e^{ik})^2 - \left( \left(V_{i,e}^R I_{ik,e}^R + V_{i,e}^I I_{ik,e}^I\right)^2 + \left(V_{i,e}^I I_{ik,e}^R - V_{i,e}^R I_{ik,e}^I\right)^2 \right) = 0 \quad (19)$$

and a corresponding set of dual variables $(\lambda_{S,e}, \mu_{S,e})$ and dual equation as (18).

## 4.5. Feasibility of the Homotopy Path

A strict condition for convergence of optimization problem with homotopy methods requires that *there exists a feasible solution (satisfying both the primal and dual feasibility) to the problem at each incremental step in the homotopy path $c(v)$ for all homotopy factors $(v \in [0,1])$.*

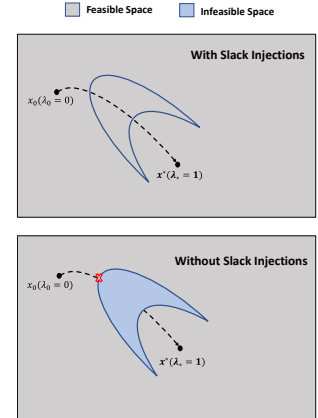

*Figure 3: Trace of homotopy curve with and without homotopy-factor dependent slack injections.*

Existing homotopy methods for power grids fail to ensure this criterion [19],[21]-[22]. To satisfy this criterion in our proposed homotopy method, we draw inspiration from the physics of the grid. In the AC-OPF problem formulation (see Section 3), the equality constraints represent the conservation of charge at each node in the grid (modeled through KCL) and the inequality constraints enforce the device limits.

Generally, it is possible that there may exist a PDIP sub-problem corresponding to a homotopy factor $v_k$ in homotopy path $(c(v_k))$, whose solution may be "infeasible" when the required current cannot be supplied by the system to satisfy KCL while ensuring that the devices are within their physical limits. If such a scenario occurs during the homotopy path $c(v_k)$, the homotopy method would fail as no feasible solution exists for a given sub-problem $k$ (see Figure 3 bottom).

To overcome this challenge, we design the IMB homotopy method to guarantee a feasible solution throughout the homotopy path by adding slack injections $(I_F^R, I_F^I)$ (also see [28]) to each node in the network that are multiplied by the homotopy factor $v$ ($vI_{F,i}^R$ and $vI_{F,i}^I$). These represent variable current sources at each node in the network to satisfy any violation of KCL. The magnitudes are minimized in the objective function with a large weight to sufficiently penalize any slack currents. The KCL equations at each node in the network given by (1b)-(1c) in the AC-OPF formulation are modified during homotopy to include the slack currents as follows:

$$\hat{I}_{G,i}^R - \hat{I}_{D,i}^R = \mathbf{Re}\left\{\sum_{k=1}^{|B|}(V_{ik}^R + jV_{ik}^I)(\hat{G}_{ik}^Y + j\hat{B}_{ik}^Y)\right\} + vI_{F,i}^R \quad \forall i \in B \quad (20)$$

$$\hat{I}_{G,i}^I - \hat{I}_{D,i}^I = \mathbf{Im}\left\{\sum_{k=1}^{|B|}(V_{ik}^R + jV_{ik}^I)(\hat{G}_{ik}^Y + j\hat{B}_{ik}^Y)\right\} + vI_{F,i}^I \quad \forall i \in B \quad (21)$$

By adding this homotopy factor based slack injections, we can vary the magnitudes of slack injections to always have a feasible solution for each step in the homotopy path. Equivalently, a feasible network implies that the Jacobian ($\mathcal{H}'(v)$) is full rank at each point on homotopy curve $c(v)$ as KCL is satisfied at each node due to the existence of slack injection sources. We further ensure that these injections are removed in the final sub-problem ($v = 0$) to solve the original network, as they are parameterized by the homotopy factor $v$. Note that this assumes there exists at least one feasible solution to the original problem $(x^*, \lambda^*, \mu^*)$ that satisfies (3a) through (3d). In case there is a problem that does not have a feasible solution (one that satisfies (1b)-(11) at $v = 0$), then one can obtain its physically infeasible solution by satisfying only relaxed form of (1b)-(1c) by replacing $vI_{F,i}^{I,R}$ in (20)-(21) with $(v + e)I_{F,i}^{I,R}$ where $e$ is a very small number.

### 4.6. Limiting for Primal-Dual Feasibility

For each homotopy step increment, we solve a PDIP sub-problem that must satisfy (3a) through (3f). In the IMB approach, we satisfy (3a)-(3d) through solving the set of non-linear equations following the damped Newton's method. We then satisfy (3e)-(3f) during each NR step through the technique described below.

Primal feasibility (3e) is ensured when the primal variables, $x$, do not violate their appropriate device limits, $h(x) \leq 0$. Similarly, the dual feasibility (3f) is maintained when the dual variables are non-negative, $\mu \geq 0$. To achieve the following, we employ a limiting heuristic similar to that of diodes in circuit simulation [27]. The limiting method ensures that any limited primal variable $x_i$ does not violate its limits (satisfies (3e)) during each NR iteration. To do so, every variable $x_i$ in the vector $x$, has an associated damping factor, $\tau_{x,i} \in (0,1]$ that modifies the updated variable, $x_i^{k+1}$ at the $k^{th}$ iteration.

$$x_i^{k+1} = x_i^k + \tau_{x,i}\Delta x_i^k \quad (22)$$

where $x_i^{k+1}$ is the value at the next iteration that is determined by the step, $\Delta x_i^k$. To properly limit the variables during each iteration to satisfy the inequality constraint ($x < \bar{x}$), the damping factor is calculated as

$$\tau_{x,i} = \min(1, \alpha_x \frac{\bar{x}_i - x_i^k}{\Delta x_i^k}) \quad (23)$$

where $\alpha_x$ is a constant (typically of 0.95-0.99) to ensure that $x_i^{k+1}$ does not hit its limit exactly. Importantly, each variable being limited calculates its own separate damping factor. Similarly, corresponding dual variables are limited to remain non-negative within each NR step:

$$\mu_i^{k+1} = \mu_i^k + \tau_{\mu,i}\Delta\mu_i^k \quad (24)$$

and vector of damping factor $(\tau_\mu)$ is calculated following this logic:

$$\tau_{\mu,i} = \min(1, -\alpha_\mu \frac{\mu_i^k}{\Delta\mu_i^k}) \quad (25)$$

### 4.7. Convergence Notes on Homotopy Methods

Global convergence of homotopy methods require following conditions to be met [16]:

i. *Defined path for the homotopy method i.e. $c(v) \in \mathcal{H}^{-1}(0)$ with $(x, v) \in range(c)$ must be smooth and should exist.*
ii. *If a curve $c$ exists, then it should intersect the final solution at $v = 0$.*

In the proposed IMB homotopy method, the first condition can be met through implicit function theorem and requires that the Jacobian ($\mathcal{H}'(v)$) of the homotopy function is of full rank for all values of $v$ along the curve. For AC-OPF, we ensure the full-rank of the iteration matrix by following the methodology in Section 4.5.

The second condition is more easily met and is linked to existence theorems in non-linear analyses [16]. If some boundary condition exists that prevents the curve from extending to infinity prior to intersecting the solution at $v = 0$, then this condition is met. In our formulation, different limiting methods ensure that the solution at any point on the curve $c$ does not diverge and extend to infinity given a feasible solution exist for the defined problem.

## 5. SUGAR Optimization Algorithm with IMB Homotopy

Using the IMB homotopy algorithm with convergence guarantees under the hold of certain conditions, we integrated the method within the circuit-based solver power flow solver, SUGAR [17] to solve the AC-OPF problem. The framework begins by accepting an input file(s) to parse the network data, generator costs and device operating limits. Following this step, all system states are initialized based on the values in the input file or based on other heuristics and any out-of-bound states are reset to feasible values within the bounds. Next, the device models are generated for elements in the network file and a homotopy parameter is embedded within the models (see, [17] and [5] for details). For non-linear models, the non-linear Hessian terms are initialized based on the initial states. The perturbed KKT conditions corresponding to the grid optimization problem at initial homotopy factor ($v$=1) are iteratively solved with Newton method, thereby providing the search direction to convergence. The homotopy factor is updated and the process is repeated with the final solution from the previous homotopy factor step used as the initial condition for the new sub-problem with updated homotopy factor. Within each NR iteration, the residual and error corresponding to the primal and dual constraints given in (3a)-(3d) are used as input heuristics to determine the homotopy factor step and the Newton step size. Also limiting heuristic based on Section 4.6 is used to ensure primal and dual feasibility given by (3e)-(3f). The algorithm is terminated once convergence is achieved for the primal and dual constraints at homotopy factor of 0 and primal and dual variables are in feasible space.

**ALGORITHM**

1: **procedure:**
2: **parse** network data, generator cost, and device limits
3: **initialize** system states to be in the feasible region; assign tolerance
4: **modify** the device models and bounds to initialize homotopy ($v$=1)
5: **create** the matrix structure
6: **loop** through all the linear devices and **create** linear admittance matrix
7: **loop** through all the non-linear devices and **initialize** Hessian matrix of non-linear terms based on initial state
  **while** homotopy factor ($v$≠0)
8:   **while** not converged **do:**
       **Solve** linearized KKT with homotopy factor, $v$, using NR
7:     **enforce** primal and dual feasibility
8:     **check** error and residual profile
9:     **if** error or residual profile diverging**:**
10:       **adjust** homotopy factor ($v$) (decrease homotopy step)
11:     **end if**
12:     **if** change in $v$: **update** homotopy elements
13:     **update** non-linear devices
14:   **end while;** update homotopy factor (increase the step)
15: **end while**

TABLE 1: AC-OPF SOLUTION COMPARISON.

| Case | # Node | SUGAR | | | fmincon [30] | | | MIPS [14] | | | IPOPT [15] | | |
|---|---|---|---|---|---|---|---|---|---|---|---|---|---|
| | | Obj. Val. ($) | Time (sec) | # Iter. | Obj. Val. ($) | Time (sec) | # Iter. | Obj. Val. ($) | Time (sec) | # Iter | Obj. Val. ($) | Time (sec) | # Iter. |
| 2869pegase[1] | 2869 | 4.58e5 | 29.6 | 243 | 4.58e5 | 13.3 | 29 | 4.58e5 | 11.2 | 41 | 4.58e5 | 13.84 | 47 |
| 2869pegase[3] | 2869 | 2.32e7 | 53.7 | 451 | 2.32e7 | 28.3 | 59 | 2.32e7 | 10.9 | 38 | 2.32e7 | 14.5 | 41 |
| 9241pegase[1] | 9241 | 1.20e6 | 135.9 | 350 | 1.20e6 | 117 | 46 | 1.20e6 | 43.8 | 51 | 1.20e6 | 40.73 | 46 |
| 9241pegase[2] | 9241 | 3.15e5 | 134.3 | 332 | 3.15e5 | 295.1 | 75 | 3.15e5 | 39.4 | 40 | 3.15e5 | 104.8 | 144 |
| 13659pegase[1] | 13659 | 1.46e6 | 299.3 | 489 | 1.46e6 | 315.1 | 65 | NC | NC | NC | NC | NC | NC |
| 13659pegase[2] | 13659 | 3.85e5 | 302.4 | 472 | 3.85e5 | 1962 | 85 | 3.85e5 | 88.8 | 79 | NC | NC | NC |
| ACTIVSg10k[1] | 10000 | 7.86e5 | 124.6 | 299 | 7.86e5 | 210 | 157 | NC | NC | NC | 7.86e5 | 127 | 182 |
| ACTIVSg10k[2] | 10000 | 1.53e5 | 122.5 | 283 | 1.53e5 | 191.9 | 108 | NC | NC | NC | 1.53e5 | 129.6 | 132 |
| ACTIVSg25k[1] | 25000 | 1.21e6 | 205.4 | 219 | 1.21e6 | 525.1 | 190 | NC | NC | NC | 1.21e6 | 198.9 | 100 |
| ACTIVSg25k[2] | 25000 | 2.39e5 | 226.2 | 241 | 2.39e5 | 307.6 | 270 | NC | NC | NC | 2.39e5 | 160.1 | 84 |

| | | | | | | | | | | | | |
|---|---|---|---|---|---|---|---|---|---|---|---|---|
| ACTIVSg70k[1] | 70000 | 3.08e6 | 1707. | 569 | 3.08e6 | 2648 | 191 | NC | NC | NC | 3.08e6 | 1377 | 264 |
| ACTIVSg70k[2] | 70000 | 6.06e5 | 1910. | 668 | 6.06e5 | 1903 | 204 | NC | NC | NC | 6.06e5 | 850.3 | 289 |
| SyntheticUSA[1] | 82000 | 4.17e6 | 1846. | 429 | NC | NC | NC | NC | NC | NC | 4.17e6 | 1702. | 236 |
| SyntheticUSA[C] | 82000 | 4.18e6 | 4610. | 865 | NC | NC | NC | NC | NC | NC | NC | NC | NC |
| ACTIVSg10k[C] | 10000 | 1.53e5 | 167.2 | 337 | NC | NC | NC | NC | NC | NC | NC | NC | NC |

1. NC refers to Non-convergent cases
2. Cases shaded in grey represent congested networks also superscripted using C
3. Superscripts 1, 2, and 3 represent three different objective functions. All the input .m files are uploaded to public git repo.

## 6. Results

To demonstrate the efficacy of our approach, we optimize the AC dispatch (with AC-OPF) for multiple test networks ranging from 9k+ nodes pegase network [30] to 80k+ node synthetic test networks for the U.S. interconnection [31]. We run two sets of analyses on each network with one corresponding to a set of non-congested scenarios (unshaded in the table) where none of the transmission line limits are binding and another corresponding to a set of congested scenarios (shaded in grey in the table) where some of the transmission line limits are binding. We run these two sets of scenarios (in total 15 cases) and compare the obtained results against those obtained from three other solvers with the Hessian supplied through a user-defined function:

- fmincon [30] solver within Matpower 7.0 [29] framework
- MIPS (Matpower Interior Point Solver) [14] solver within Matpower 7.0 [29] framework
- IPOPT (Interior Point Optimizer) [15] solver within Matpower 7.0 [29] framework

The following changes to the default settings are made for the commercial solvers: i) Increase the maximum iteration count for the IPOPT solver from 250 to 750; ii) Represent maximum flow limits in terms of currents (I) instead of default setting of real power (P) for congested cases; and iii) For non-congested scenarios, use default settings for flow limits of real power (this is because convergence was much worse when flow limits in terms of currents were used).

All of the .m testcase files with generation cost data, line flow limits and other case information are made publicly available in a Gitlab repository [29]. Also, all the output .raw files with updated states from SUGAR-IMB based solver are added to the same repository.

The condensed results with comparison against standard non-linear optimization tools are shown in Table 1. The SUGAR-based approach with the IMB homotopy method demonstrates better convergence across all (15) scenarios compared to fmincon that converged for 12 out of 15 scenarios, IPOPT that converged for 11 out of 15 scenarios and MIPS that converged for 5 out of 15 scenarios. While the SUGAR-IMB homotopy solver does take more iterations compared to other tools to reach an optimal solution, the net runtime for the SUGAR-IMB solver is lower in many of those large scenarios when compared against other methods. This is primarily because the homotopy methods used in SUGAR-IMB ensure smooth convergence for each sub-problem. The standard commercial tools, however, perform more internal function evaluations between each iteration to converge to the final solution. Also, it's worth noting that SUGAR-IMB outperforms other methods for congested large networks (as these tend to have more complex solution space) as shown in the bottom two rows of the results table.

## 7. Conclusions

In this paper, we develop a methodology for robustly solving hard-to-solve AC-OPF problems especially when good initial conditions are not available. The proposed framework introduces a homotopy method called Incremental Model Building approach that initially modifies the problem to solve a trivial problem and then gradually constrains the problem to solve the AC-OPF exactly in the final homotopy factor. The homotopy method has its roots in circuit-theory and with inclusion of grid physics can satisfy critical criteria for convergence of homotopy methods. To demonstrate the robustness of the proposed framework, we show results for hard-to-solve networks and compare those against other standard non-convex optimization engines. The Incremental Model Building framework acts as a basis for future planning optimizations of the grid, such as optimal power flow and transmission expansion, that require a robust engine capable of scaling to large systems regardless of initial conditions.

## 8. Acknowledgement

This work was supported in part by the National Science Foundation under grant ECCS-1800812.